\documentclass[11pt,fleqn]{article}
\usepackage{amsmath,amssymb,graphicx}

\begin{document}
\sloppy
\newtheorem{axiom}{Axiom}[section]
\newtheorem{conjecture}[axiom]{Conjecture}
\newtheorem{corollary}[axiom]{Corollary}
\newtheorem{definition}[axiom]{Definition}
\newtheorem{example}[axiom]{Example}
\newtheorem{fact}[axiom]{Fact}
\newtheorem{lemma}[axiom]{Lemma}
\newtheorem{observation}[axiom]{Observation}
\newtheorem{open}[axiom]{Problem}
\newtheorem{proposition}[axiom]{Proposition}
\newtheorem{theorem}[axiom]{Theorem}

\renewcommand{\topfraction}{1.0}
\renewcommand{\bottomfraction}{1.0}

\newcommand{\proof}{\emph{Proof.}\ \ }
\newcommand{\qed}{~~$\Box$}
\newcommand{\rz}{{\mathbb{R}}}
\newcommand{\nz}{{\mathbb{N}}}
\newcommand{\zz}{{\mathbb{Z}}}
\newcommand{\eps}{\varepsilon}
\newcommand{\cei}[1]{\lceil #1\rceil}
\newcommand{\flo}[1]{\left\lfloor #1\right\rfloor}

\newcommand{\boxxx}[1]
 {\begin{center}\fbox{\begin{minipage}{10.00cm}#1\smallskip\end{minipage}}\end{center}}

\newcommand{\aaa}{\alpha}
\newcommand{\bbb}{\beta}
\newcommand{\ccc}{\gamma}
\newcommand{\ddd}{\delta}
\newcommand{\jjj}{J^{\pi}}
\newcommand{\xxx}{L}
\newcommand{\yyy}{K}
\newcommand{\lam}{\lambda}
\newcommand{\qqq}{\gamma}
\newcommand{\gmam}{GMAM}

\title{{\bf Well-solvable cases of the QAP with block-structured matrices}}
\author{
\sc Eranda \c{C}ela\thanks{{\tt cela@opt.math.tu-graz.ac.at}.
Institut f\"ur Optimierung und Diskrete Mathematik, TU Graz, Steyrergasse 30, A-8010 Graz, Austria}
\and\sc Vladimir G.\ Deineko\thanks{{\tt Vladimir.Deineko@wbs.ac.uk}.
Warwick Business School, The University of Warwick, Coventry CV4 7AL, United Kingdom}
\and\sc Gerhard J.\ Woeginger\thanks{{\tt gwoegi@win.tue.nl}.
Department of Mathematics and Computer Science, TU Eindhoven, P.O.\ Box 513,
5600 MB Eindhoven, Netherlands}
}
\date{}
\maketitle

\begin{abstract}
We investigate special cases of the quadratic assignment problem (QAP) where one of the two underlying
matrices carries a simple block structure.
For the special case where the second underlying matrix is a monotone anti-Monge matrix, we derive a
polynomial time result for a certain class of cut problems.
For the special case where the second underlying matrix is a product matrix, we identify two sets of
conditions on the block structure that make this QAP polynomially solvable respectively NP-hard.

\end{abstract}

\medskip
\section{Introduction}
The \emph{Quadratic Assignment Problem} (QAP) is an important and well-studied problem in
combinatorial optimization; we refer the reader to the book \cite{Cela-book} by \c{C}ela
and the recent book \cite{Burkard-book} by Burkard, Dell'Amico \& Martello for comprehensive
surveys on this problem.
The QAP in Koopmans-Beckmann form \cite{KoBe1957} takes as input two $n\times n$ square matrices
$A=(a_{ij})$ and $B=(b_{ij})$ with real entries.
The goal is to find a permutation $\pi$ that minimizes the objective function
\begin{equation}
\label{eq:qap}
Z_{\pi}(A,B) ~:=~ \sum_{i=1}^n\sum_{j=1}^n~ a_{\pi(i)\pi(j)} \, b_{ij}.
\end{equation}
Here $\pi$ ranges over the set $S_n$ of all permutations of $\{1,2,\ldots,n\}$.
In general, the QAP is extremely difficult to solve and hard to approximate.
One branch of research on the QAP concentrates on the algorithmic behavior of strongly structured
special cases; see for instance Burkard \& al \cite{BCRW1998}, Deineko \& Woeginger \cite{DeWo1998},
\c{C}ela \& al \cite{Cela2011}, or \c{C}ela, Deineko \& Woeginger \cite{Cela2012} for typical results
in this direction.
We will contribute several new results to this research branch.

\paragraph{Results of this paper.}
Our first result is motivated by a \emph{balanced} multi-cut problem where a group of entities 
has to be divided into $q$ clusters of equal size with the objective of minimizing the overall 
connection cost between different clusters; we refer the reader to Lengauer \cite{Lengauer-book} 
for a variety of applications of this and other cut problems in VLSI design.
Pferschy, Rudolf \& Woeginger \cite{PfRuWo1994} show that if the connection costs carry a certain
anti-Monge structure, the balanced multi-cut problem can be solved efficiently.
In the language of the QAP, the connection costs can be summarized in an anti-Monge matrix $A$ and
the balanced multi-cut structure can be encoded by a block-structured $0$-$1$ matrix $B$.
The result of \cite{PfRuWo1994} then states that the corresponding special case of the QAP is
always solved by the identity permutation.

In this paper, we generalize the result of \cite{PfRuWo1994} to the case where the entities are to
be divided into $q$ clusters of prescribed (but not necessarily equal) sizes.
In our generalization the connection costs are given by a \emph{monotone} anti-Monge matrix $A$,
and the cluster structure is specified by a block-structured $0$-$1$ multi-cut matrix $B$ that 
lists the clusters in order of non-decreasing size.
We show that the resulting special case of the QAP (again) is always solved by the identity permutation.
Our proof method strongly hinges on the monotonicity of matrix $A$; in fact it can be seen that without 
monotonicity the result would break down (see Example~\ref{ex:non-monotone} in Section~\ref{sec:cut}).

Our second result concerns a wide class of specially structured QAPs that are loosely related to
the multi-cut QAP in the preceding paragraph.
Matrix $A$ is now a product matrix, and hence a special anti-Monge matrix with a particularly nice
structure.
Matrix $B$ is now a block matrix with some fixed block pattern $P$, and hence a generalization of
the multi-cut matrices in the multi-cut QAP.
In comparison to the multi-cut QAP, the structure of matrix $A$ has become more restricted, while 
the structure of matrix $B$ has become more general.
The resulting special case of the QAP is called the \emph{Product-Block QAP} with block pattern $P$.

On the positive side, we identify conditions on the block pattern $P$ that render the Product-Block
QAP polynomially solvable.
One main ingredient of the polynomial algorithm is the concavity of certain underlying functions,
and the other main ingredient is an extensive enumeration of cases.
On the negative side, we identify conditions on the block pattern $P$ that make the Product-Block
QAP NP-hard.
The positive conditions as well as the negative conditions on the pattern exploit the connections to
an underlying continuous quadratic program.

\paragraph{Organization of the paper.}
Section~\ref{sec:definitions} introduces all the relevant matrix classes, and also states some
simple observations on the QAP.
Section~\ref{sec:cut} contains our results on the multi-cut problem on anti-Monge matrices, and
Section~\ref{sec:problo} presents our results on the Product-Block QAP.
Section~\ref{sec:conclusions} concludes the paper by listing some open problems.

\medskip
\section{Definitions and preliminaries}
\label{sec:definitions}
All matrices in this paper are symmetric and have real entries.
In order to avoid trouble with the standard models of computation, we will sometimes assume for 
our complexity results that the matrix entries are \emph{rational} numbers; this assumption
will always be stressed and stated explicitly in the corresponding theorem.

For a $q\times q$ matrix $P=(p_{ij})$, we say that an $n\times n$ matrix $B=(b_{ij})$ is a
\emph{block matrix with block pattern $P$} if the following holds:
(i) there exists a partition of the row and column set $\{1,\ldots,n\}$ into $q$ (possibly empty)
intervals $I_1,\ldots,I_q$ such that for $1\le k\le q-1$ all elements of interval $I_k$ are smaller 
than all elements of interval $I_{k+1}$;
(ii) for all indices $i$ and $j$ with $1\le i,j\le n$ and $i\in I_k$ and $j\in I_{\ell}$, we have
$b_{ij}=p_{k\ell}$.
The sets $I_1,\ldots,I_q$ form the so-called row and column blocks of matrix $B$.

A \emph{multi-cut matrix} $B$ is a block matrix whose pattern matrix has $0$'s along the main 
diagonal and $1$'s everywhere else.
Intuitively speaking, every block $I_k$ in a multi-cut matrix represents a cluster of data points;
data points in the same cluster are very similar to each other (and hence at distance $0$), whereas
data points from different clusters are dissimilar and far away from each other.
A multi-cut matrix is in \emph{normal form}, if its block sizes are in non-decreasing order with
$|I_1|\le|I_2|\le\cdots\le|I_q|$; note that the rows and columns of every multi-cut matrix can 
easily be permuted into this normal form.

For a real number $\lam>0$, a \emph{1-$\lam$-1 block matrix} is a block matrix with the following
block pattern $P(\lam)$:
\begin{equation}
\label{eq:rays}
P(\lam) ~=~ \left( \begin{array}{ccc} 0&0&0\\0&0&1\\0&1&\lam \end{array} \right)
\end{equation}
Note that for $0<\mu<\lam$, any 1-$\lam$-1 block matrix can be written as a non-negative linear
combination of two 1-$\mu$-1 block matrices.

An $n\times n$ matrix $A=(a_{ij})$ is \emph{monotone}, if $a_{ij}\le a_{i,j+1}$ and
$a_{ij}\le a_{i+1,j}$ holds for all $i,j$, that is, if the entries in every row and every column
are in non-decreasing order.
Matrix $A$ is a \emph{sum matrix}, if there are (not necessarily positive) real numbers
$\aaa_1,\ldots,\aaa_n$ such that $a_{ij}=\aaa_i+\aaa_j$ for $1\le i,j\le n$.
Matrix $A$ is a \emph{product matrix}, if there are \emph{non-negative} real numbers
$\aaa_1,\ldots,\aaa_n$ such that
\begin{equation}
\label{eq:product}
a_{ij} ~=~ \aaa_i\,\aaa_j \mbox{\qquad for $1\le i,j\le n$.}
\end{equation}
If $\aaa_1\le\aaa_2\le\cdots\le\aaa_n$ holds, then (\ref{eq:product}) gives a \emph{monotone} product matrix.
Matrix $A$ is an \emph{anti-Monge matrix}, if its entries are non-negative and satisfy the anti-Monge inequalities
\begin{equation}
\label{eq:anti-Monge}
a_{ij}+a_{rs} ~\ge~ a_{is}+a_{rj} \mbox{\qquad for $1\le i<r\le n$ and $1\le j<s\le n$.}
\end{equation}
In other words, in every $2\times2$ submatrix the sum of the entries on the main diagonal dominates
the sum of the entries on the other diagonal.
This property essentially dates back to the work of Gaspard Monge \cite{Monge} in the 18th century.
Much research has been done on the effects of Monge structures in combinatorial optimization, and we
refer the reader to the survey \cite{BuKlRu1996} by Burkard, Klinz \& Rudolf for more information on
Monge and anti-Monge structures.
It can be shown (see Rudolf \& Woeginger \cite{RuWo1995}; Burkard \& al \cite{BCRW1998}) that a
symmetric matrix is a \emph{monotone} anti-Monge matrix, if and only if it can be written as a
non-negative linear combination of 1-2-1 block matrices; in other words, the 1-2-1 block matrices
in (\ref{eq:rays}) form the extremal rays of the cone of monotone anti-Monge matrices.
Furthermore every (arbitrary, not necessarily monotone) anti-Monge matrix can be written as a sum of
a monotone anti-Monge matrix and an appropriately chosen sum matrix.

A matrix belongs to the \emph{$\lam$-generalized monotone anti-Monge cone}, $\lam$-{\gmam} cone for short,
if it can be written as a non-negative linear combination of 1-$\lam$-1 block matrices.
Note that the $2$-{\gmam} cone coincides with the standard monotone anti-Monge cone, and that for
$0<\mu<\lam$ the $\mu$-{\gmam} cone properly contains the $\lam$-{\gmam} cone.

\bigskip
We close this section with two simple but useful results from the QAP folklore.
\begin{fact}
\label{fact:folklore.1}
Consider a QAP where $A$ is a sum matrix and where all row and column sums of $B$ are the same.
Then all permutations yield the same objective value.
\end{fact}
\proof
Assume that the entries of matrix $A$ are given by $a_{ij}=\aaa_i+\aaa_j$ for $1\le i,j\le n$,
and let $\bbb$ denote the row and column sum of matrix $B$.
Then for every permutation $\pi$, the objective value in (\ref{eq:qap}) equals
$Z_{\pi}(A,B)=2\bbb\sum_{i=1}^n\aaa_i$.
\qed

\begin{fact}
\label{fact:folklore.2}
If a permutation $\pi$ solves the QAP between matrices $A'$ and $B$ and the QAP between matrices
$A''$ and $B$ to optimality, then $\pi$ also solves the QAP between $A'+A''$ and $B$ to optimality.
\qed
\end{fact}

\medskip
\section{Multi-cut problems on anti-Monge matrices}
\label{sec:cut}
In this section, we consider the special case of the QAP where matrix $A$ is a monotone anti-Monge
matrix and where matrix $B$ is a multi-cut matrix in normal form.
We will first resolve a highly restricted special case in Section~\ref{sec:cut.1}, and then
deduce a polynomial time result for the general version in Section~\ref{sec:cut.2}.

\subsection{A highly restricted special case}
\label{sec:cut.1}
In this section we consider the special case of the QAP where the $n\times n$ matrix $A$ is a 1-2-1
block matrix and where the $n\times n$ matrix $B$ is a multi-cut matrix with only two blocks 
(and hence a standard cut matrix).
We denote the sizes of the three blocks of matrix $A$ by $r,s,t$, and the sizes of the two blocks of
matrix $B$ by $u$ and $v$; see Figure~\ref{fig:cut} for an illustration.
As $B$ is in normal form, we have $u\le v$.
Note furthermore that $r+s+t=u+v=n$, and that $s+t\le2v$.

\begin{figure}[bt]
\begin{center}
\unitlength=0.8mm
\begin{picture}(60,57)(-5,-1)
\thicklines
\put( 0, 0){\line(0,1){50}}
\put(20, 0){\line(0,1){50}}
\put(40, 0){\line(0,1){50}}
\put(50, 0){\line(0,1){50}}
\put( 0, 0){\line(1,0){50}}
\put( 0,10){\line(1,0){50}}
\put( 0,30){\line(1,0){50}}
\put( 0,50){\line(1,0){50}}
\put(10, 5){\makebox(0,0)[cc]{\large 0}}
\put(30, 5){\makebox(0,0)[cc]{\large 1}}
\put(45, 5){\makebox(0,0)[cc]{\large 2}}
\put(10,20){\makebox(0,0)[cc]{\large 0}}
\put(30,20){\makebox(0,0)[cc]{\large 0}}
\put(45,20){\makebox(0,0)[cc]{\large 1}}
\put(10,40){\makebox(0,0)[cc]{\large 0}}
\put(30,40){\makebox(0,0)[cc]{\large 0}}
\put(45,40){\makebox(0,0)[cc]{\large 0}}
\put(-4, 5){\makebox(0,0)[cc]{\large $t$}}
\put(-4,20){\makebox(0,0)[cc]{\large $s$}}
\put(-4,40){\makebox(0,0)[cc]{\large $r$}}
\put(10,53){\makebox(0,0)[cc]{\large $r$}}
\put(30,53){\makebox(0,0)[cc]{\large $s$}}
\put(45,53){\makebox(0,0)[cc]{\large $t$}}
\end{picture}
\qquad\qquad
\begin{picture}(60,57)(-5,-1)
\thicklines
\put( 0, 0){\line(0,1){50}}
\put(20, 0){\line(0,1){50}}
\put(50, 0){\line(0,1){50}}
\put( 0, 0){\line(1,0){50}}
\put( 0,30){\line(1,0){50}}
\put( 0,50){\line(1,0){50}}
\put(10,15){\makebox(0,0)[cc]{\large 1}}
\put(35,15){\makebox(0,0)[cc]{\large 0}}
\put(10,40){\makebox(0,0)[cc]{\large 0}}
\put(35,40){\makebox(0,0)[cc]{\large 1}}
\put(-4,15){\makebox(0,0)[cc]{\large $v$}}
\put(-4,40){\makebox(0,0)[cc]{\large $u$}}
\put(10,53){\makebox(0,0)[cc]{\large $u$}}
\put(35,53){\makebox(0,0)[cc]{\large $v$}}
\end{picture}
\caption{The notation $r,s,t,u,v$ for the 1-2-1 block matrix $A$ and the multi-cut matrix $B$
used in Section~\protect{\ref{sec:cut.1}}.}
\label{fig:cut}
\end{center}
\end{figure}

Consider an arbitrary permutation $\pi$ for the formulation (\ref{eq:qap}) of the QAP, and let $x$
(respectively, $y$ and $z$) be the number of rows from $A$'s first block (respectively, its second and
third block) that $\pi$ assigns to the first block of matrix $B$; the remaining $r-x$ (respectively,
$s-y$ and $t-z$) rows from these blocks are assigned to the second block of matrix $B$.
The corresponding objective value can then be written as
\begin{eqnarray}
f(x,y,z) &=& 2y(t-z) +2z(s-y) +4z(t-z) \nonumber \\[0.8ex]
&=& -4z^2-4yz+2ty+(2s+4t)z.  \label{eq:woe.1}
\end{eqnarray}
The variables $x,y,z$ are integers with $0\le x\le r$, $0\le y\le s$, and $0\le z\le t$, that satisfy $x+y+z=u$.
Of course, the optimal objective value of the QAP coincides with the optimal objective value of this
non-linear integer program IP.

Since the value of $f(x,y,z)$ in (\ref{eq:woe.1}) does not depend on $x$, we will drop variable $x$ from
our further considerations.
Since $x=u-y-z$ and $u+v=r+s+t$, the constraint $0\le x\le r$ can be rewritten as $s+t-v\le y+z\le u$.
Next we define a continuous programming relaxation CPR of the IP, in which $y$ and $z$ are real variables.
Furthermore, we relax the upper bound constraint $y+z\le u$ to the less restrictive constraint $y+z\le v$.
Therefore the CPR has the objective of minimizing (\ref{eq:woe.1}) subject to the constraints
\begin{equation}
\label{eq:woe.2}
0\le y\le s; \mbox{\qquad}
0\le z\le t; \mbox{\qquad}
s+t-v\le y+z\le v.
\end{equation}
The following two lemmas derive lower bounds on the optimal objective value of the CPR, and hence also
on the optimal objective value of the IP and the QAP.

\begin{lemma}
\label{le:woe.1}
If $t\le v\le s+t$, then the objective value of the CPR (and hence also of the QAP) is at
least $\qqq_1:=2t\,(s+t-v)$.
\end{lemma}
\proof
We start with two auxiliary inequalities.
The condition $t\le v$ yields
\begin{equation}
\label{eq:woe.1a}
2ts ~\ge~ 2t\,(s+t-v) ~=~ \qqq_1,
\end{equation}
and $s+t\le2v$ yields
\begin{equation}
\label{eq:woe.1b}
2(2v-s)\,(s+t-v) ~\ge~ 2t\,(s+t-v) ~=~ \qqq_1.
\end{equation}
Since the Hessian matrix of function $f$ is indefinite, its minimizers lie on the boundary of the
feasible region defined by (\ref{eq:woe.2}).
We distinguish six cases on the six bounding lines.

(Case~1).
The minimizer satisfies $y=0$.
The problem turns into the minimization of $g(z)=-4z^2+(2s+4t)z=2z(s+2t-2z)$ subject to $0\le z\le t$ and
$s+t-v\le z\le v$, and hence subject to $s+t-v\le z\le t$.
As $g(z)$ is concave, it is minimized at the boundary.
Inequalities (\ref{eq:woe.1a}) and (\ref{eq:woe.1b}) show that $g(t)=2st$ and $g(s+t-v)=2(2v-s)\,(s+t-v)$
are both at least $\qqq_1$.

(Case~2).
The minimizer satisfies $y=s$.
The problem is to minimize $g(z)=2(t-z)(2z+s)$ subject to $0\le z\le t$ and $t-v\le z\le v-s$, and hence
subject to $0\le z\le v-s$.
Again $g(z)$ is concave, and (\ref{eq:woe.1a}) and (\ref{eq:woe.1b}) show that the values $g(0)=2st$ and
$g(v-s)=2(s+t-v)(2v-1)$ at the boundary are at least $\qqq_1$.

(Case~3).
The minimizer satisfies $z=0$.
The problem is to minimize $g(y)=2ty$ subject to $0\le y\le s$ and $s+t-v\le y\le v$.
Then $y\ge s+t-v$ implies $g(y)\ge2t\,(s+t-v)=\qqq_1$.

(Case~4).
The minimizer satisfies $z=t$.
The problem is to minimize $g(y)=2t(s-y)$ subject to $0\le y\le s$ and $s-v\le y\le v-t$.
Then $y\le v-t$ implies $g(y)\ge2t\,(s+t-v)=\qqq_1$.

(Case~5).
The minimizer satisfies $y+z=s+t-v$.
The problem is to minimize the increasing linear function $g(z)=2(2v-s-t)z+2(s+t-v)t$ subject to
$0\le z\le t$ and $t-v\le z\le s+t-v$.
Then $z\ge0$ yields $g(z)\ge2(s+t-v)t=\qqq_1$.

(Case~6).
The minimizer satisfies $y+z=v$.
The problem is to minimize the decreasing linear function $g(z)=2vt-2(2v-s-t)z$ subject to
$0\le z\le t$ and $v-s\le z\le v$.
Then $z\le t$ yields $g(z)\ge2t(s+t-v)=\qqq_1$.
\qed

\begin{lemma}
\label{le:woe.2}
If $v<t$, then the objective value of the CPR (and hence also of the QAP) is at
least $\qqq_2:=2v\,(s+2t-2v)$.
\end{lemma}
\proof
This proof is analogous to the proof of the preceding lemma, except that we use a different set
of bounds and inequalities.
{From} $s+t\le2v$ we conclude
\begin{equation}
\label{eq:woe.2a}
2(s+t-v)(2v-s) ~=~ s\,(2v-s-t) +\qqq_2 ~\ge~ \qqq_2.
\end{equation}
This time the feasible region is bounded by only four straight lines.
The case $z=0$ is impossible, since then $s+t-v\le y+z=y\le s$ implies the contradiction $t\le v$.
Also the case $z=t$ is impossible, since then $t=z\le y+z\le v$ implies $t\le v$.
Hence we distinguish only four cases on the four bounding lines.

(Case~1).
The minimizer satisfies $y=0$.
The problem is to minimize $g(z)=-4z^2+(2s+4t)z=2z(s+2t-2z)$ subject to $0\le z\le t$ and $s+t-v\le z\le v$,
and hence subject to $s+t-v\le z\le v$.
Function $g(z)$ is concave, and we have $g(v)=2v(s+2t-2v)=\qqq_2$, and (\ref{eq:woe.2a}) implies that
$g(s+t-v)=2(s+t-v)(2v-s)\ge\qqq_2$.

(Case~2).
The minimizer satisfies $y=s$.
The problem is to minimize $g(z)=2(t-z)(2z+s)$ subject to $0\le z\le t$ and $t-v\le z\le v-s$, and hence
subject to $t-v\le z\le v-s$.
Once again $g(z)$ is concave, and $g(t-v)=2v(s+2t-2z)=\qqq_2$, and (\ref{eq:woe.2a}) implies that
$g(v-s)=2(s+t-v)(2v-s)\ge\qqq_2$.

(Case~3).
The minimizer satisfies $y+z=s+t-v$.
The problem is to minimize the increasing linear function $g(z)=2(2v-s-t)z+2(s+t-v)t$ subject to
$0\le z\le t$ and $t-v\le z\le s+t-v$.
Then $z\ge t-v$ yields $g(z)\ge g(t-v)=\qqq_2$.

(Case~4).
The minimizer satisfies $y+z=v$.
The problem is to minimize the decreasing linear function $g(z)=2vt-2(2v-s-t)z$ subject to
$0\le z\le t$ and $v-s\le z\le v$.
Then $z\le v$ yields $g(z)\ge g(v)=\qqq_2$.
\qed

\begin{theorem}
\label{th:woe.rays}
If $A$ is a 1-2-1 block matrix and $B$ is a multi-cut matrix in normal form with two blocks,
then the identity permutation solves the QAP to optimality.
\end{theorem}
\proof
If $v\ge s+t$, then the identity permutation yields an objective value of $0$ which clearly is optimal.
If $t\le v\le s+t$, then the identity permutation sets $y=s+t-v$ and $z=0$ and yields objective value $\qqq_1$;
by the lower bound in Lemma~\ref{le:woe.1} this is optimal.
If $v<t$, then the identity permutation sets $y=0$ and $z=t-v$ and yields objective value $\qqq_2$;
by the lower bound in Lemma~\ref{le:woe.2} this is optimal.
\qed

\subsection{The general case}
\label{sec:cut.2}
Now we are ready to establish our main result for the general case, where the multi-cut 
matrix $B$ has an arbitrary number of blocks.

\begin{theorem}
\label{th:woe.main}
If $A$ is a monotone anti-Monge matrix and $B$ is a multi-cut matrix in normal form,
then the identity permutation solves the QAP to optimality.
\end{theorem}
\proof
The proof is done in two steps.
In the first step, we assume that matrix $A$ is a 1-2-1 block matrix.
Let $I_1,\ldots,I_q$ with $|I_1|\le|I_2|\le\cdots\le|I_q|$ denote the blocks of matrix $B$.
Consider an optimal permutation $\pi^*$ for the QAP, and assume that $\pi^*$ assigns the row
set $J_k$ of $A$ to block $I_k$ of $B$, where $1\le k\le q$.
The submatrix $A'$ of $A$ induced by the rows and columns in $J_k\cup J_{k+1}$ and the submatrix $B'$
induced by the rows and columns in $I_k\cup I_{k+1}$ satisfy the conditions of Theorem~\ref{th:woe.rays}.
According to the theorem, we may repartition $J_k$ and $J_{k+1}$ such that all elements in $J_k$ precede
all the elements in $J_{k+1}$, as imposed by the identity permutation.
Repeated application of such repartitioning eventually transforms $\pi^*$ into the identity permutation
without worsening the objective value.

In the second step, we consider the most general case with an arbitrary monotone anti-Monge matrix $A$.
As $A$ can be written as a non-negative linear combination of 1-2-1 block matrices, and as the identity
permutation optimally solves the QAP between any 1-2-1 block matrix and matrix $B$, the identity
permutation also optimally solves the QAP between $A$ and $B$ according to Fact~\ref{fact:folklore.2}.
\qed

\bigskip
Next, we want to demonstrate that Theorem~\ref{th:woe.main} generalizes the following proposition from
the Monge literature.
\begin{proposition}
\label{pr:PRW}
(Pferschy, Rudolf \& Woeginger \cite{PfRuWo1994})\\
If $A$ is a symmetric (not necessarily monotone) anti-Monge matrix and $B$ is a multi-cut matrix
with all blocks of identical size, the identity permutation solves the QAP to optimality.
\end{proposition}
\proof
First note that all row and column sums in matrix $B$ are identical.
The anti-Monge matrix $A$ can be written as the sum of a monotone anti-Monge matrix $A'$ and a sum matrix $A''$.
Theorem~\ref{th:woe.main} yields that the identity permutation optimally solves the QAP between $A'$ and $B$,
and Fact~\ref{fact:folklore.1} yields that every permutation (and in particular the identity permutation)
optimally solves the QAP between $A''$ and $B$.
\qed

\bigskip
The following two examples illustrate that in a certain sense the statement in Theorem~\ref{th:woe.main}
is best possible.
If we allow $A$ to be a general and not necessarily monotone anti-Monge matrix, the statement fails.
If we take $A$ from a slightly larger $\lam$-generalized monotone anti-Monge cone $\lam$-{\gmam} with 
$\lam<2$ (and do not restrict it to the standard cone with $\lam=2$), the statement fails.

\begin{example}
\label{ex:non-monotone}
Consider the QAP with the following non-monotone anti-Monge matrix $A$ and the following 
multi-cut matrix $B$ in normal form:
$$
A=\left( \begin{array}{ccc} 2&1&1\\1&0&0\\1&0&0 \end{array} \right) \mbox{\qquad\qquad}
B=\left( \begin{array}{ccc} 0&1&1\\1&0&0\\1&0&0 \end{array} \right)
$$
Then the identity permutation has objective value $4$, whereas the permutation that switches the
first row (and column) of $A$ with its third row (and column) has a better objective value of $2$.
\end{example}

\begin{example}
\label{ex:other-cones}
For $\lam<2$, consider the QAP with the following matrix $A$ in the $\lam$-{\gmam} cone and the
following multi-cut matrix $B$ in normal form:
$$
A=\left( \begin{array}{cccc} 0&0&1&1\\0&0&1&1\\1&1&\lam&\lam\\1&1&\lam&\lam \end{array} \right)
\mbox{\qquad\qquad}
B=\left( \begin{array}{cccc} 0&0&1&1\\0&0&1&1\\1&1&0&0\\1&1&0&0 \end{array} \right)
$$
Then the identity permutation has objective value $8$, whereas the permutation that switches the
second row (and column) of $A$ with its fourth row (and column) has a better objective value of $4+2\lam$.
\end{example}

Finally, we show that another natural extension of the special case in Theorem~\ref{th:woe.main}
yields an NP-hard problem.
\begin{lemma}
\label{le:woe.np}
The QAP remains NP-hard,
even if $A$ is a monotone matrix and $B$ is a multi-cut matrix.
\end{lemma}
\proof
The proof is done by means of a reduction from the {\sc Graph Bisection} problem
(see for instance Garey \& Johnson \cite{GaJo1979}) which is known to be NP-hard.
The input for {\sc Graph Bisection} consists of an undirected graph $G=(V,E)$ on $n=2r$ vertices,
together with an integer bound $t$.
The goal is to decide whether there exists a partition of $V$ into two subsets $V_1$ and $V_2$
with $|V_1|=|V_2|=r$, such that at most $t$ edges in $E$ connect $V_1$ to $V_2$.

We construct the following QAP instance.
The $n\times n$ matrix $A$ is the sum of the adjacency matrix $A'$ of graph $G$ and of the sum matrix $A''$
that is defined by $a''_{ij}=2i+2j$ for $1\le i,j\le n$; note that $A$ indeed is monotone.
The $n\times n$ matrix $B$ is the multi-cut matrix with two blocks of size $r$.
It is straightforward to verify with the help of Facts~\ref{fact:folklore.1} and~\ref{fact:folklore.2}
that the {\sc Graph Bisection} instance has answer YES, if and only if the constructed QAP instance has
a feasible solution with objective value at most $n^2(n+1)+t$.
\qed

\medskip
\section{The Product-Block QAP}
\label{sec:problo}
In this section we study the so-called \emph{Product-Block QAP}, the special case where matrix $A$
is a product matrix and where matrix $B$ is a block matrix with some fixed pattern $P$.
Throughout this section we assume that all patterns (and hence all considered block matrices) have rational entries.
We stress that pattern $P$, and in particular the number $q$ of blocks in $P$, are not part of the input.
(Note that if $P$ is part of the input, then one may choose $P=B$ so that matrix $B$ essentially
remains unrestricted.)
We remind the reader that matrix $P$ and all other matrices in this paper are symmetric.

The following definitions play a central role in our investigations.
A \emph{bad ensemble} for a $q\times q$ pattern matrix $P=(p_{ij})$ consists of the following:
\begin{itemize}
\itemsep=-0.1ex
\item two indices $r$ and $s$ with $1\le r<s\le q$,
\item a real number $\ccc$ with $0\le\ccc\le1$,
\item real numbers $\ell_i$ with $0\le\ell_i\le1$ for $i\in\{1,\ldots,q\}\setminus\{r,s\}$.
\end{itemize}
With every bad ensemble, we associate the quadratic program QP-1 for non-negative real variables
$x_1,\ldots,x_q$ in Figure~\ref{fig:QP-1}.
We will only consider ensembles for which the feasible region specified by
(\ref{eq:QP-1.1})--(\ref{eq:QP-1.3}) is non-empty.
The crucial property of a bad ensemble is that QP-1 has a unique minimizer $(x^*_1,\ldots,x^*_q)$, 
and that this minimizer satisfies
\begin{equation}
\label{eq:bad}
0<x^*_r<\ccc \mbox{\qquad and\qquad} 0<x^*_s<\ccc.
\end{equation}
In a similar spirit, we introduce \emph{very bad ensembles} for a $q\times q$ pattern matrix $P=(p_{ij})$
that consist of the following:
\begin{itemize}
\itemsep=-0.1ex
\item two indices $r$ and $s$ with $1\le r<s\le q$,
\item real numbers $\ell_i$ and $u_i$ with $0\le\ell_i<u_i\le1$ for $i=1,\ldots,q$.
\end{itemize}
With a very bad ensemble, we associate the quadratic program QP-2 in Figure~\ref{fig:QP-2}.
We stress that the upper bound constraints $x_i<u_i$ in QP-2 are strict.
The crucial properties of a very bad ensemble are that QP-2 has a \emph{unique} minimizer 
$(x^*_1,\ldots,x^*_q)$, that all the $x^*_i$ are rational, and that 
\begin{equation}
\label{eq:verybad}
\ell_r<x^*_r<u_r \mbox{\qquad and\qquad} \ell_s<x^*_s<u_s.
\end{equation}

\begin{figure}[bt]
\boxxx{\begin{subequations}
\begin{alignat}{4}
\min~~ &\sum_{i=1}^q\sum_{j=1}^q\,p_{ij}\,x_ix_j \nonumber \\[1.0ex]
s.t.~~ &\sum_{i=1}^qx_i =1   \label{eq:QP-1.1}     \\[1.0ex]
       &x_r+x_s = \ccc       \label{eq:QP-1.2}     \\[1.0ex]
       &x_i=\ell_i \mbox{\quad\qquad for $i\in\{1,\ldots,q\}\setminus\{r,s\}$} \label{eq:QP-1.3}
\\[-6.0ex] \nonumber
\end{alignat}
\end{subequations}}
\vspace{-2ex}
\caption{The continuous quadratic program QP-1.}
\label{fig:QP-1}
\end{figure}

\begin{figure}[bt]
\vspace{4ex}
\boxxx{\begin{subequations}
\begin{alignat}{4}
\min~~ &\sum_{i=1}^q\sum_{j=1}^q\,p_{ij}\,x_ix_j \nonumber \\[1.0ex]
s.t.~~ &\sum_{i=1}^qx_i =1   \label{eq:QP-2.1}     \\[1.0ex]
       &\ell_i\le x_i<u_i \mbox{\qquad for $i=1,\ldots,q$.} \label{eq:QP-2.2}
\\[-6.0ex] \nonumber
\end{alignat}
\end{subequations}}
\vspace{-2ex}
\caption{The continuous quadratic program QP-2.}
\label{fig:QP-2}
\vspace{1ex}
\end{figure}

In Sections~\ref{sec:problo.1} and~\ref{sec:problo.2} we will prove the following theorem, which
settles the computational complexity of the Product-Block QAP for two large families of patterns.
\begin{theorem}
\label{le:pro.main}
Consider the Product-Block QAP with a fixed rational block pattern $P$.
\begin{itemize}
\itemsep=-0.1ex
\item[(i)]  If $P$ does not allow any bad ensemble, the QAP is polynomially solvable.
\item[(ii)] If $P$ has a very bad ensemble, the QAP is NP-hard.
\end{itemize}
\end{theorem}

In general, it is not straightforward to see whether a given pattern matrix allows a bad or very bad ensemble.
The following two corollaries extract two clean and tidy pattern classes that are covered by
Theorem~\ref{le:pro.main}.

\begin{corollary}
\label{co:pro.poly}
The Product-Block QAP is polynomially solvable whenever the rational pattern matrix $P=(p_{ij})$ satisfies
\begin{equation}
\label{eq:light}
p_{ii}+p_{jj} ~\le~ 2p_{ij} \mbox{\qquad for $1\le i,j\le n$.}
\end{equation}
\end{corollary}
\proof
Assume for the sake of contradiction that there is a bad ensemble and consider the corresponding
minimizer $(x^*_1,\ldots,x^*_q)$.
As every variable $x_i$ with $i\notin\{r,s\}$ is frozen at $x_i=\ell_i$ and as variable $x_s$ can be replaced
by $\ccc-x_r$ according to (\ref{eq:QP-1.2}), QP-1 boils down to minimizing the quadratic function
\[ g(x_r) ~=~ (p_{rr}+p_{ss}-2p_{rs})\,x_r^2 +c_1\,x_r + c_0 \]
subject to the constraint $0\le x_r\le\ccc$; here the coefficients $c_0$ and $c_1$ are certain real
numbers that depend on the ensemble.
Condition (\ref{eq:light}) implies that the coefficient of the quadratic term $x_r^2$ is non-positive,
so that $g(x_r)$ is concave and takes its minimum at $x_r=0$ or $x_r=\ccc$.
This contradicts (\ref{eq:bad}), and hence Theorem~\ref{le:pro.main}.(i) applies.
\qed

\begin{corollary}
\label{co:pro.hard}
The Product-Block QAP is NP-hard whenever there exist two indices $r$ and $s$, for which
the rational pattern matrix $P=(p_{ij})$ satisfies
\begin{equation}
\label{eq:heavy}
p_{rr}>p_{rs} \mbox{\quad and\quad} p_{ss}>p_{rs}.
\end{equation}
\end{corollary}
\proof
We construct a very bad ensemble for the $2\times2$ submatrix of $P$ spanned by rows (and columns) $r$ and $s$.
If we set $\ell_r=\ell_s=0$ and $u_r=u_s=1$, routine calculations show that QP-2 is minimized at
the rational point
\begin{equation}
\label{eq:heavy.1}
x^*_r ~=~  \frac{p_{ss}-p_{rs}}{p_{rr}+p_{ss}-2p_{rs}}
\mbox{\quad and\quad}
x^*_s ~=~  \frac{p_{rr}-p_{rs}}{p_{rr}+p_{ss}-2p_{rs}}.
\end{equation}
By (\ref{eq:heavy}) the numerators and denominators in (\ref{eq:heavy.1}) are positive,
so that $0<x^*_r,x^*_s<1$.
Hence this ensemble indeed is very bad and Theorem~\ref{le:pro.main}.(ii) applies.
\qed

\bigskip
As the quadratic programs QP-1 and QP-2 are easy to analyze for $q=2$, some routine calculations
yield the following corollary.
\begin{corollary}
\label{co:pro.2x2}
The Product-Block QAP with a $2\times2$ rational block pattern $P=(p_{ij})$ is NP-hard 
if $p_{11}>p_{12}$ and $p_{22}>p_{12}$, and it is polynomially solvable otherwise.
\qed
\end{corollary}

\subsection{Proof of the polynomial time result}
\label{sec:problo.1}
In this section we prove the positive statement in Theorem~\ref{le:pro.main}.(i).
Throughout we consider an $n\times n$ product matrix $A=(\aaa_i\,\aaa_j)$; we assume without loss of
generality that $\aaa_1\le\aaa_2\le\cdots\le\aaa_n$ so that $A$ is monotone, and we furthermore assume
that the values $\aaa_i$ are normalized so that $\sum_{i=1}^n\aaa_i=1$.
For a subset $J$ of $\{1,\ldots,n\}$ we denote $\aaa(J)=\sum_{i\in J}\aaa_i$.
Matrix $B$ is a block matrix with a $q\times q$ pattern $P$ that does not allow any bad ensemble;
the blocks of $B$ are $I_1,\ldots,I_q$.

Consider a permutation $\pi$ that for $k=1,\ldots,q$ assigns the row set $\jjj_k$ of $A$ to block $I_k$ of $B$.
We say that the sets $\jjj_i$ and $\jjj_j$ are \emph{separable}, if either all elements of $\jjj_i$ are less
than or equal to all elements of $\jjj_j$, or if all elements of $\jjj_i$ are greater than or equal to all
elements of $\jjj_j$.

\begin{lemma}
\label{le:pro.separable}
There exists an optimal permutation $\pi$, such that for all $r\ne s$ the sets $\jjj_r$ and $\jjj_s$ are separable.
\end{lemma}
\proof
Let $\pi$ be an optimal permutation.
We denote $y_k=\aaa(\jjj_k)$ for $1\le k\le q$, and we observe that the objective value of the QAP can
be written as $\sum_{i=1}^q\sum_{j=1}^q p_{ij}y_iy_j$.
We define an ensemble with $\ccc=y_r+y_s$, and with $\ell_i=u_i=y_i$ for $i\notin\{r,s\}$.
In QP-1 we freeze every variable $x_i$ with $i\notin\{r,s\}$ at its current value $x_i=y_i$, and we
furthermore substitute $x_s=\ccc-x_r$.
The resulting quadratic program asks to minimize a uni-variate quadratic function $g(x_r)$ subject to
the single constraint $0\le x_r\le\ccc$.
As the pattern $P$ does not allow bad ensembles, the minimum is taken at the boundary.
Hence function $g(x_r)$ is either concave over $[0,\ccc]$, or it is convex and its minimizer lies
outside of $[0,\ccc]$.

For the QAP this means that the objective value is minimized if $\aaa(\jjj_r)$ either becomes as small
or as large as possible.
As matrix $A$ is monotone, this in turn means that set $\jjj_r$ should consist either of the $|\jjj_r|$
smallest or the $|\jjj_r|$ largest elements in $\jjj_r\cup\jjj_s$, and that set $\jjj_s$ should consist of
the remaining elements.
In other words, we are able to separate the sets $\jjj_r$ and $\jjj_s$ without worsening the objective value.
Repeated application of this separation step will eventually transform permutation $\pi$ into the
desired form.
\qed

\begin{theorem}
\label{th:pro.algorithm}
If the $q\times q$ pattern $P$ does not allow any bad ensemble, the Product-Block QAP with block
pattern $P$ is solvable in time $O(q^2q!+n\log n)$.
\end{theorem}
\proof
By Lemma~\ref{le:pro.separable}, there are only $q!$ many cases to check.
After sorting the numbers $\aaa_i$ and after performing some appropriate preprocessing, the objective
value for every such case can be determined in $O(q^2)$ time.
\qed

\bigskip
As $q$ is not part of the input, Theorem~\ref{th:pro.algorithm} completes the proof of
Theorem~\ref{le:pro.main}.(i).
Note that Theorem~\ref{th:woe.main} shows that a \emph{single} permutation is optimal for all instances
in the considered multi-cut problem, whereas Theorem~\ref{th:pro.algorithm} leaves us with a huge number 
of $q!$ candidate permutations for the Product-Block QAP.
The following example illustrates that even for a $2\times2$ pattern $P$, there is no way of sharpening
our statements to a single optimal permutation.

\begin{example}
Matrices $A_1$ and $A_2$ are monotone product matrices, and Corollary~\ref{co:pro.poly}
shows that the pattern of block matrix $B$ does not allow any bad ensemble:
$$
A_1=\left( \begin{array}{ccc} 1&1&2\\1&1&2\\2&2&4 \end{array} \right) \mbox{\qquad}
A_2=\left( \begin{array}{ccc} 1&2&2\\2&4&4\\2&4&4 \end{array} \right) \mbox{\qquad}
  B=\left( \begin{array}{ccc} 0&2&2\\2&1&1\\2&1&1 \end{array} \right)
$$
By Lemma~\ref{le:pro.separable}, there are two candidates for the optimal permutation: the identity
permutation $\pi_1$ and the permutation $\pi_2$ that assigns the first/second/third row (and column) 
of matrix $A_i$ respectively to the third/second/first row (and column) of matrix $B$.
For the QAP between $A_1$ and $B$, permutation $\pi_1$ with objective value $21$ loses to permutation
$\pi_2$ with objective value $20$.
But for the QAP between $A_2$ and $B$, permutation $\pi_1$ with objective value $32$ beats permutation
$\pi_2$ with objective value $33$.
\end{example}

\subsection{Proof of the hardness result}
\label{sec:problo.2}
In this section we prove the negative statement in Theorem~\ref{le:pro.main}.(ii).
Throughout we consider a fixed $q\times q$ pattern $P$ with a very bad ensemble, and its corresponding
unique, rational minimizer $(x^*_1,\ldots,x^*_q)$.
The optimal objective value of the quadratic program QP-2 is denoted
$z^*=\sum_{i=1}^q\sum_{j=1}^q\,p_{ij}\,x^*_ix^*_j$.

If $x^*_i=0$ holds for some $i$, then we may ignore this variable and its corresponding row and column
in the pattern in our further considerations; hence we will assume without loss of generality that
$x^*_i>0$ holds for all $i$.
We fix a large integer $\yyy$ such that $x^*_i\yyy$ is integer for all $i$, and such that 
\begin{equation}
\label{eq:hard.1a}
\yyy ~>~ \max\left\{\frac{2}{u_i-x^*_i},\, \frac{1}{x^*_i}\right\} \mbox{\qquad for $i=1,\ldots,q$.}
\end{equation}
and
\begin{equation}
\label{eq:hard.1b}
\yyy ~>~ \frac{1}{x^*_j-\ell_j}  \mbox{\qquad\quad for $j\in\{r,s\}$.}
\end{equation}

The NP-hardness proof is done by means of a reduction from the following variant of the number 
partition problem (see Garey \& Johnson \cite{GaJo1979}).
An instance of {\sc Partition} consists of a sequence $v_1,\ldots,v_m$ of positive rational numbers
with $\sum_{k=1}^mv_k=2$.
The goal is to decide whether there exists a subset $M\subset\{1,\ldots,m\}$ such that
$\sum_{k\in M}v_k=1$ and $\sum_{k\notin M}v_k=1$.

We begin the reduction by introducing the following numbers.
We define an integer $\xxx=m\yyy$. 
Our choice of $\yyy$ in (\ref{eq:hard.1a}) yields $K>1/x^*_i$, and consequently
\begin{equation}
\label{eq:hard.2}
x^*_i\xxx ~>~ m \mbox{\qquad\quad for all $i$.}
\end{equation}
The bound (\ref{eq:hard.1a}) together with $\xxx\ge\yyy$ implies $x^*_i+2/\xxx<u_i$ for all $i$, and
the bound (\ref{eq:hard.1b}) with $\xxx\ge\yyy$ implies $x^*_j-1/\xxx>\ell_j$ for $j\in\{r,s\}$.

Now let us construct an instance of the Product-Block QAP.
The dimension is chosen as $n=\xxx-2$.
The product matrix $A$ is specified by the following $\aaa$-values:
\begin{itemize}
\itemsep=-0.1ex
\item For $k=1,\ldots,m$, we introduce a partition-value $\aaa_k=(1+v_k)/\xxx$.
\item For $k=m+1,\ldots,n$, we introduce a dummy-value $\aaa_k=1/\xxx$.
\end{itemize}
Note that every $\aaa$-value is at least~$1/\xxx$, and that the overall sum of all the $\aaa$-values 
equals $(n+2)/\xxx=1$.
The block matrix $B$ has blocks $I_1,\ldots,I_q$ and obeys the block pattern $P$.
Block $I_r$ has size $x^*_r\xxx-1$, block $I_s$ has size $x^*_s\xxx-1$, and every remaining block 
$I_i$ with $i\notin\{r,s\}$ has size $x^*_i\xxx$.
This completes the description of the QAP instance.

\begin{lemma}
\label{le:hard.1}
If the {\sc Partition} instance has answer YES, then for the constructed instance of the QAP
has a permutation with objective value at most $z^*$.
\end{lemma}
\proof
Let $M\subset\{1,\ldots,m\}$ be a solution for the {\sc Partition} instance.
We define a permutation $\pi$ that assigns the $\aaa$-values of matrix $A$ to the blocks of matrix $B$:
\begin{itemize}
\item
To block $I_r$, we assign the $|M|$ partition-values $\aaa_k$ with $k\in M$ together with 
$x^*_r\xxx-|M|-1$ dummy-values.
Note that the number $x^*_r\xxx-|M|-1$ is non-negative by (\ref{eq:hard.2}), and note that 
the sum of all assigned $\aaa$-values is $x^*_r$.
\item
To block $I_s$, we assign the remaining $m-|M|$ partition-values $\aaa_k$ with $k\notin M$ together
with $x^*_s\xxx-m+|M|-1$ dummy-values.
Note that $x^*_s\xxx-m+|M|-1$ is non-negative by (\ref{eq:hard.2}), and note that the sum of all 
assigned $\aaa$-values is $x^*_s$.
\item
To every remaining block $I_i$ with $i\notin\{r,s\}$, we assign $x^*_i\xxx$ dummy-values.
The sum of all assigned $\aaa$-values is $x^*_i$.
\end{itemize}
The resulting objective value for the QAP is $\sum_{i=1}^q\sum_{j=1}^q\,p_{ij}\,x^*_ix^*_j$, and
hence coincides with $z^*$.
\qed

\begin{lemma}
\label{le:hard.2}
If the constructed instance of the QAP has a permutation with objective value at most $z^*$,
then the {\sc Partition} instance has answer YES.
\end{lemma}
\proof
Consider a permutation $\pi$ for the QAP with objective value at most $z^*$.
Let $y_i$ denote the sum of all $\aaa$-values that $\pi$ assigns to block $I_i$, and
let $\ddd_i=y_i-|I_i|/\xxx$.
It is easily seen that $\ddd_i\ge0$ for all $i$ and that $\sum_{i=1}^q\ddd_i=2/\xxx$, which implies
$0\le\ddd_i\le2/\xxx$.
For $i\notin\{r,s\}$ we have $|I_i|=x^*_i\xxx$, and hence $x^*_i\le y_i\le x^*_i+(2/\xxx)$.
For $j\in\{r,s\}$ we have $|I_j|=x^*_j\xxx-1$, and hence $x^*_j-(1/\xxx)\le y_j\le x^*_j+(1/\xxx)$.
By our choice of $\yyy$ and $\xxx$, these bounds imply $\ell_i\le y_i<u_i$ for all $i$.

This means that $(y_1,\ldots,y_q)$ constitutes a feasible solution for the quadratic program QP-2 
with objective value $\sum_{i=1}^q\sum_{j=1}^qp_{ij}y_iy_j\le z^*$, and hence is a minimizer for QP-2.
As we are working with a very bad ensemble, the minimizer is unique so that $y_i=x^*_i$ for all~$i$.
This leads to $\ddd_i=0$ for $i\notin\{r,s\}$, and $\ddd_r=\ddd_s=1/\xxx$.
Consequently permutation $\pi$ assigns all partition-values $\aaa_k$ to the two blocks $I_r$ and $I_s$.
If we define set $M$ to contain all indices $k$ for which the partition-value $\aaa_k$ is assigned to 
block $I_r$, it is easily seen that $\sum_{k\in M}v_k=1$.
Hence the {\sc Partition} instance has answer YES.
\qed

\bigskip
Lemma~\ref{le:hard.1} and Lemma~\ref{le:hard.2} together establish the correctness of our reduction.
This completes the proof of Theorem~\ref{le:pro.main}.(ii).

\medskip
\section{Conclusions}
\label{sec:conclusions}
We have studied a family of special cases of the quadratic assignment problem, where one
matrix carries an anti-Monge structure and where the other matrix has a simple block structure.
We identified a number of well-behaved cases that are solvable in polynomial time,
and we also got some partial understanding of the borderline between certain easy and hard cases.
Many questions remain open, and we will now list some of them.

Examples~\ref{ex:non-monotone} and~\ref{ex:other-cones} describe scenarios that can not be 
settled with the methodology of Section~\ref{sec:cut}.
The precise complexity of these two scenarios remains open:

\begin{open}
Consider the QAP where $A$ is a (general, not necessarily monotone) anti-Monge
matrix and where $B$ is a multi-cut matrix.
Is this special case polynomially solvable?
\end{open}

\begin{open}
For some fixed $0<\lam<2$, consider the QAP where $A$ lies in the generalized
monotone anti-Monge cone $\lam$-{\gmam} and where $B$ is a multi-cut matrix.
Is this special case NP-hard?
\end{open}

Our results on the Product-Block QAP in Section~\ref{sec:problo} leave a considerable
gap between the polynomially solvable area and the NP-hard area.
This gap can be narrowed somewhat by allowing non-rational minimizers in the proof of
Theorem~\ref{le:pro.main}.(ii), and by working with sufficiently precise rational 
approximations of all the involved real numbers; the technical details, however, would be gory. 
Corollaries~\ref{co:pro.poly} and~\ref{co:pro.hard} might indicate some vague connection to 
totally positive and totally negative matrices (see for instance Pinkus \cite{Pinkus-book}) 
and to Eigenvalue spectra.  

It is not clear that the gap actually can be closed, as for some patterns the problem 
might neither be polynomially solvable nor NP-hard.
Schaefer's famous dichotomy theorem \cite{Schaefer1978} states that every member of a large 
family of constraint satisfaction problems is either polynomially solvable or NP-complete.
Is there a similar result for the Product-Block QAP?
\begin{open}
Is there a dichotomy theorem for the Product-Block QAP showing that every pattern $P$
gives rise to either a polynomially solvable or an NP-complete problem?
\end{open}

\bigskip
\small
\paragraph{Acknowledgements.}
Part of this research was conducted while Vladimir Deineko and Gerhard Woeginger were visiting TU Graz, and were supported by the Austrian Science Fund (FWF): W1230, Doctoral Program ``Discrete Mathematics''.
Vladimir Deineko acknowledges support
by Warwick University's Centre for Discrete Mathematics and Its Applications (DIMAP).
Gerhard Woeginger acknowledges support
by DIAMANT (a mathematics cluster of the Netherlands Organization for Scientific Research NWO),
and by the Alexander von Humboldt Foundation, Bonn, Germany.
\normalsize


\end{document}